\documentclass{CSML}

\def\dOi{13(1:4)2017}
\lmcsheading%
{\dOi}
{1--6}
{}
{}
{Jun.~\phantom03, 2016}
{Feb.~\phantom02, 2017}
{}

\usepackage{amsfonts}
\usepackage{amsmath}
\usepackage{amsthm}
\usepackage{color}
\usepackage{eufrak}
\usepackage{hyperref}
\usepackage{graphicx}

\newcommand\np{\mbox{NP}}
\newcommand\pp{\mbox{P}}
\newcommand\sz{\mbox{SIZE}}

\newcommand\ut{1^{(t)}}
\newcommand\un{1^{(n)}}
\newcommand\um{1^{(n^{1/3k})}}

\newcommand\upk{{\mbox{UP}_{k}}}
\newcommand\upkconstant{{\mbox{UP}_{k,c}}}

\newcommand\upkconstanttilde{{\text{\emph{UP}}_{k,\tilde c}}}

\newcommand{\eqdef}{\stackrel{\rm def}{=}}

\def\colorful{1}

\ifnum\colorful=1

\newcommand{\blue}[1]{{{#1}}}
\newcommand{\red}[1]{{{#1}}}

\fi
\ifnum\colorful=0

\newcommand{\blue}[1]{{{#1}}}
\newcommand{\red}[1]{{{#1}}}

\fi

\begin{document}

\title[Unprovability of circuit upper bounds\\in Cook's theory PV]
      {Unprovability of circuit upper bounds\\in Cook's theory PV}

\author[J.~Kraj\'{\i}\v{c}ek]{Jan Kraj\'{\i}\v{c}ek}
\address{Faculty of Mathematics and Physics,
Charles University in Prague}
\email{\{krajicek, igor.oliveira\}@karlin.mff.cuni.cz}
\author[I.~C.~Oliveira]{Igor C.~Oliveira}
\address{\vspace{-18 pt}}

\thanks{Supported in part by CNPq grant 200252/2015-1.}

\begin{abstract}
We establish unconditionally
that for every integer $k \geq 1$ there is a language $L \in \pp$
such that it is consistent with Cook's theory PV that $L \notin \sz(n^k)$. 
Our argument is non-constructive and does not provide an explicit description of this language.
\end{abstract}

\maketitle
\section{Introduction}
Bounded arithmetic theories constitute a class of weak subtheories of
Peano arithmetic with close ties to computational complexity theory.
Prominent among them is theory PV defined by Cook \cite{Coo75}
as an equational theory and later reformulated as a universal
first order theory in \cite{KPT,kniha}.

Theory PV or its mild extensions seem to formalize most 
of contemporary complexity theory (cf. 
\cite{PW,Bus-book,kniha,Jer04,Jer07,Jer09,CooNgu,Pich-phd,Pich,Pich-pcp}
and references therein). \red{For instance, it is known that the PCP Theorem can be formalized and proved in PV \cite{Pich-pcp}.} 
It is thus of interest to understand, given an established
conjecture, whether it is provable in one of these theories or
at least consistent with them. 

An unprovability statement 
can be understood as a result illustrating the
inadequacy of methods available in the 
respective theory. This is studied in complexity theory as the so
called barriers \red{(cf. \cite{BGS75,RR97,AW09})}, often formulated using ad hoc concepts
hard to compare with each other. The unprovability results
on the other hand
are in the tradition how mathematical logic captured (and answered)
similar questions in other parts of mathematics. 

The latter direction, to show the consistency of the conjecture in question
with PV or with stronger theories, is at least as interesting as 
showing its unprovability.
Such a consistency result says that, although we do not know if the conjecture
is true (meaning true in the standard model of natural numbers), we
know that it is true in a non-standard model of a theory so strong that
complexity theory looks in it almost indistinguishable from the standard one.

In this work we study the provability of circuit upper bounds (or equivalently,
the consistency of lower bounds). Circuit lower bounds were considered 
in bounded arithmetic by 
Razborov \cite{Raz95} in a particular formalism. We use the somewhat more
intrinsic formalism of \cite{PW,Bus-book,kniha} and followed in 
\cite{Jer04,CooKra,CooNgu,Pich-phd,Pich,Pich-pcp}.

It has been proved in \cite{CooKra},
assuming that $\np \not\subseteq$~co$\np/O(1)$ or 
that the polynomial time hierarchy does not collapse to the Boolean hierarchy,
that it is consistent with PV that 
$\np \not\subseteq  \pp/$poly.
Here we prove {\em unconditionally} that for every $k \geq 1$ 
there is a language $L \in \pp$
such that it is consistent with Cook's theory PV that $L \notin \sz(n^k)$, where $\sz(n^k)$ denotes
the class of languages decided by \emph{non-uniform} Boolean circuits of size at most \red{$O(n^k)$}. We refer to the statement of Theorem \ref{t:main} below for a precise formulation of the result.

We do not know how to extend our result to Buss's theory $S^1_2$ from
\cite{Bus-book} (results from \cite{CooKra} were extended there to $S^1_2$)
or how to show that one can take SAT for $L$ for all $k \geq 1$.
Perhaps the most accessible problem is to extend our result to 
PV augmented by the dual weak pigeonhole principle for polynomial time functions,
a theory denoted APC$_1$ by some authors.

\section{Formalization and statement of the theorem}\label{s:statement}

The language of PV has function symbols for all polynomial time
algorithms as generated by Cobham's limited recursion on notation \cite{Cob64}.
All axioms of PV are universal formulas codifying how particular
algorithms are defined from each other. The details of the definition of
PV are fairly technical, but such details are needed only for establishing links
between PV and propositional proof systems (cf. \cite{kniha}). 
We use a form of Herbrand's Theorem
(see below), and for that it only matters that the axioms are universal formulas.
In fact, we could add to PV any set of true universal sentences as additional axioms,
and our unprovability result would still hold.

We will talk about polynomial time algorithms in the theory meaning that they
are represented by the corresponding function symbols. We shall claim on a few 
occasions that some algorithm $f_1$ constructed in a particular way from
another algorithm $f_2$ can be defined in PV; this means that PV proves that
$f_1$ behaves as described in the definition. In all cases this is straightforward
but tedious, and presupposes a certain amount of bootstrapping of PV
 which is part of standard
background in bounded arithmetic (see e.g.$\;$in \cite{Bus-book} how this is done).
The details are not necessary for 
understanding our argument and can be found in \cite{Bus-book,Coo75,kniha,KPT}.

For a unary PV function symbol $f$ and integers \red{$k, c \geq 1$}, denote by
\red{$\upkconstant(f)$} the sentence
\begin{equation}\label{6a}
\forall 1^{(n)} \exists \mbox{circuit }C_n (|C_n| \le \red{c}n^k) \forall x (|x|=n),\  
f(x)\neq 0 \leftrightarrow C_n(x)=1\ ,
\end{equation}
which asserts that the (polynomial time) language defined by $f$ admits a (non-uniform) sequence of circuits of size at most $\red{c} n^k$.\footnote{For the reader familiar with bounded arithmetic, we stress that we abuse notation and use $|C_n|$ to denote the number of gates in $C_n$, while $|x|$ refers to the length of $x$ in the usual sense. Also, the symbol
$\forall 1^{(n)}$ abbreviates the universal quantification over strings of the form 
$1^{(n)}$, i.e., strings consisting of a sequence of ones.} (We refer to \cite{kniha, Pich-phd} for 
more information about the formalization of circuit complexity in bounded arithmetic.)

\begin{thm}\label{t:main}
For every $k \geq 1$ there is a unary PV function symbol $h$ such that
\red{for no constant $c \geq 1$ PV proves the sentence \emph{$\upkconstant(h)$}}. 
\end{thm}

The high level idea of the proof is: (\emph{i}) the provability of (\ref{6a})
implies a certain uniformity of the family of circuits, and (\emph{ii})
we can adapt the proof by Santhanam and Williams from \cite{SW} that $\pp$ has no uniform sequences of
circuits of size \red{$O(n^k)$}, for any fixed $k \geq 1$. Complications arise
as the uniformity given by (\emph{i}) is more general than the one employed in
(\emph{ii}). In particular, it is not clear how to establish Theorem \ref{t:main} using only the soundness of PV and 
 (extensions of) the Santhanam-Williams Theorem.

 To get around this difficulty we argue roughly as follows. Either a candidate sentence $\upkconstant(g)$
that we start with is not provable in PV (and we are done), or we extract from any proof of this sentence a finite number 
of languages in $\pp$ such that PV cannot prove that all of them admit circuits of size \red{$O(n^k)$}. 
We remark that the non-constructive aspect of the result comes from the fact that the hard language and 
its deterministic time complexity may depend on a (possibly non-existent) proof of the initial sentence.

In order to implement this approach we use that PV is a universal theory for polynomial time computations, a formalization of the main ideas employed in the uniform circuit lower bound from \cite{SW}, the KPT Theorem from bounded arithmetic
(Theorem \ref{t:KPT} below), 
and a finite number of recursive applications of Herbrand's Theorem. The argument has a few subtle points, and we make some additional observations after we present the proof of Theorem \ref{t:main} in Section \ref{s:proof}.

\begin{rem}
An alternative and equally natural formalization of circuit upper bounds can be obtained via a single formula \emph{$\upk(h)$} that existentially quantifies over parameter $c$. This leads however to a sentence of higher quantifier complexity. While
KPT witnessing \emph{(}stated as Theorem \emph{\ref{t:KPT}} in Section \emph{\ref{s:proof}}\emph{)} can be generalized in this direction, the information
it then offers does not seem to yield the polynomial
time algorithms our technique needs.
In particular, we leave the unprovability of the modified version as an open problem.
\end{rem}

\red{Theorem \ref{t:main} and a standard compactness argument imply the following result.

\begin{cor}\label{c:model}
For every $k\geq 1$ there exists a unary PV function symbol $h$ and a model $\mathfrak{M}_k$ of PV such that for every $c\geq1$ we have
$$
\mathfrak{M}_k\,\models\, \neg \text{\emph{UP}}_{k,c}(h).
$$
\end{cor}

In other words, from the point of view of $\mathfrak{M}_k$ there are languages in P that require non-uniform circuits of size $\omega(n^k)$.}

\section{Uniform sequences of circuits and PV}

In this section we adapt a proof by 
Santhanam and Williams \cite{SW} that $\pp$ is not included in
$(\pp$-$\mbox{uniform})$-$\sz(n^k)$.
Here
$(\pp$-$\mbox{uniform})$-$\sz(n^k)$ is the class of languages recognizable by
a polynomial time uniform family of circuits of size at most \red{$O(n^k)$}. That is, there
is a polynomial time algorithm $f$ that on input $\un$ computes
a description of a size \red{$cn^k$} circuit $C_n$\red{, where $c \geq 1$ is a fixed constant}. Following \cite{SW}
we take as the description the set of all 4-tuples
\begin{equation}\label{1}
(\un, u, v, w)
\end{equation}
where $u, v$ are names of nodes \blue{($\le k (\log n + O(1))$ bits each)} such that there is
a wire from $u$ to $v$, and $w$ is the information about the type of the gate at $v$
or about the input at $v$ if $v$ is an input node ($\le \log n + O(1)$ bits). We assume that 
a special tuple indicates the output node of $C_n$. The language consisting of all 4-tuples (\ref{1}) for all $n \geq 1$ is 
called $L_{\text{dc}}$, the direct connection language of $\{C_n\}_n$.

The following standard definitions play an important role in the argument. 
We use $\mbox{DTIME}(n^d)/n^{2/3}$ to denote the class of languages recognizable by
a time \red{$O(n^d)$} algorithm with an advice of size $n^{2/3}$. We say that a language $L$ 
is infinitely often in a complexity class $\Gamma$ if $L$ agrees on infinitely many 
input lengths with some language $L' \in \Gamma$. 

The next lemma formalizes the deterministic time hierarchy theorem with a bounded amount of advice.

\begin{lem} \label{5}
For every $d \geq 1$ there is $L \in \mbox{DTIME}(n^{d+1})$, represented by 
algorithm $g_{d+1}$ computing its characteristic function, such that for every
time \blue{$O(n^d)$} algorithm $h$ working with advice $n^{2/3}$ there is \red{$c_h \geq 1$}
such that $\mbox{PV}$ proves\emph{:}
$$
\forall n \geq \red{c_h} \forall a (|a|= n^{2/3}) \exists x (|x|=n),\ 
h(x,a) \neq g_{d+1}(x)\ . 
$$
\end{lem}

\proof 
The separation is reported as a folklore result in \cite[Proposition
2.1]{SW}.
We simply check that its proof formalizes in PV.

Define a time \blue{$O(n^{d+1})$} algorithm $g_{d+1}$ that operates as follows. On an
input $x$ of length $|x|=n$: 
\begin{itemize}

\item it interprets the first $\log n$ bits of $x$ as a description of
a time \blue{$n^d \log n$} algorithm $h$, and the next $n^{2/3}$ bits as advice $a$,

\item runs $h$ on $x$ with advice $a$,

\item outputs $0$ if and only if the simulation ends with a non-zero value.

\end{itemize}
The constant \blue{$c_h \geq 8$} is chosen so that \red{$\log c_h$} bits suffice to describe
the particular $h$. \red{Observe that in order for the sentence to hold for \blue{every} large enough $n$ it is important that the parts of the input corresponding to the description of the algorithm and the advice are disjoint.}
\qed

Take now $\{C_n\}_n$ a $\pp$-uniform sequence of size \red{$cn^k$} circuits
and let $f$ be the generating polynomial time algorithm. That is,
on input $\un$ $f$ produces the list of 4-tuples as in (\ref{1}). Following
\cite{SW} we compress each such 4-tuple into the 5-tuple
\begin{equation}\label{2}
(\mbox{Bin}(n)0\um, u,v,w, \ut) 
\end{equation}
where $\mbox{Bin}(n)$ is the dyadic numeral for $n$ (of length $\log n + O(1)$ bits) and $t$
is chosen to pad the length of the 5-tuple to exactly $m(n) \eqdef \blue{\lceil n^{1/(2k)} \rceil}$ bits, \blue{as soon as $n$ is sufficiently large}
(parameter $t$ is not present in \cite{SW}). The language of all such 5-tuples obtained from $L_{\text{dc}}$ is the language 
$L_{\text{succ}}$, the succinct version of $L_{\text{dc}}$. It is polynomial time and an algorithm $\tilde f$ recognizing it
can be easily defined from $f$ and, in particular, in PV.

Let $\text{CircuitVal}(y,x)$ be the polynomial time algorithm
evaluating circuit $y$ on input $x$.

\begin{lem}\label{6}
Let $f$, $\{C_n\}_n$, and $\tilde f$ be as above, and assume \red{that for some $\tilde{c} \geq 1$}
\begin{equation}\label{3}
\text{\emph{PV}}\; \vdash\; \red{\upkconstanttilde(\tilde f)}\ .
\end{equation}
Let $g \eqdef g_{3k}$ for a fixed integer $k \geq 3$ be the function guaranteed to exist by Lemma \emph{\ref{5}}.
Then there exists $\red{c_f} \geq 1$ such that PV proves
\begin{equation} \label{4}
\forall 1^{(n)} (n \geq \red{c_f}) \exists x (|x|=n),\ 
g(x) \neq C_n(x)\ ,
\end{equation}
where $C_n(x)$ abbreviates $\text{\emph{CircuitVal}}(f(\un),x)$.
\end{lem}

\proof Our argument will follow the proof of \cite[Theorem 1.1]{SW} and is done in PV.
Assuming (\ref{4}) fails we describe an explicit polynomial time algorithm $h$
that will certify that $g$ is (infinitely often) in $\text{DTIME}(n^{3k-1})/n^{2/3}$. 
This contradicts the sentence from Lemma \ref{5}.

Algorithm $h$ operates as follows.  
By the assumption (\ref{3}) there are circuits $D_m$ recognizing $L_{\text{succ}}$ on $m$-bit inputs, where $m = m(n)$, as defined above.
Upon receiving $x$, $|x|=n$, and advice string 
$a$, $|a|= n^{2/3}$, describing \blue{a candidate} circuit $D_m$, 
$h$ tries all possible 3-tuples $(u,v,w)$ (among no more than $O(n^{2k+1})$ possibilities)
and for each of them uses $D_m$ to check if the corresponding 5-tuple as in
(\ref{2}) is in $L_{\text{succ}}$. Since \red{for large enough $n$ the corresponding circuit} $D_m$ has size $\blue{O(n^{1/2})}$,
this requires time $O(n)$ for each 5-tuple. There are $O(n^{2k+1})$ such 
simulations so the total time this part takes is $O(n^{2k+2})$.

After this stage $h$ knows the description of $C_n$, a circuit of size at most $cn^k$, and uses it to compute \blue{a candidate value for} $g(x)$ in time $O(n^{2k})$. \blue{Under our initial assumption,} the algorithm is correct on infinitely many input lengths, which is contradictory if $k \geq 3$.
\qed

\begin{lem}\label{7}
Let $f$, $g$, $k$, $\{C_n\}_n$, and $\tilde f$ be as above, and assume that
\emph{(\ref{3})} holds.
There is $\red{c_f} \geq 1$ and a polynomial time algorithm $e$ such that PV proves
\begin{equation}
\forall 1^{(n)} (n \geq \red{c_f}),\ 
|e(\un)|=n \wedge g(e(\un)) \neq C_n(e(\un))\ .
\end{equation}
That is, $e$ provably produces witnesses to \emph{(\ref{4})}.
\end{lem}

\proof This follows from Lemma \ref{6} using Herbrand's Theorem, as (\ref{4}) is a
$\forall\exists$-formula and PV is a universal theory.
\qed

\section{Proof of Theorem \ref{t:main}}\label{s:proof}

We will need the following standard witnessing result from bounded arithmetic
(the so called KPT theorem), stated below for convenience of the reader.

\begin{thm}[\cite{KPT}, see also \cite{kniha}]\label{t:KPT}
Let $T$ be a universal theory with vocabulary $\mathcal{L}$, $\phi$ be an open $\mathcal{L}$-formula, and suppose that
$$
T\, \vdash\, \forall w\,\exists u\,\forall v\, \phi(w,u,v)\ .
$$
Then there exist a constant $k \geq 1$ and a finite sequence $t_1, \ldots, t_k$ of $\mathcal{L}$-terms such that
$$
T\, \vdash\, \phi(w,t_1(w), v_1) \vee \phi(w,t_2(w,v_1),v_2) \vee \ldots \vee \phi(w, t_k(w, v_1, \ldots, v_{k-1}), v_k)\ ,
$$
where the notation $t_i(w, v_1, \ldots, v_{i-1})$ indicates that these are the only variables occurring in $t_i$.   
\end{thm}
We remark that Theorem \ref{t:KPT} has a natural interpretation as an interactive game with finitely many rounds, 
and we refer to \cite{Pich} for an example in the related context of circuit lower bounds.

Continuing with the proof of Theorem \ref{t:main}, assume 
\begin{equation}\label{6b}
\mbox{PV}\ \vdash\ \red{\upkconstant(g)}\ ,
\end{equation}
where $g = g_{3k}$ and \red{$c \geq 1$ is arbitrary}.  Observe that \red{$\upkconstant(\cdot)$} is a sentence of the form $\forall \exists \forall \phi$, where $\phi$ is an open formula in the language of PV. 
By Theorem \ref{t:KPT} there are polynomial time algorithms
$f_1, \dots, f_r$ where $r$ is a fixed constant such that PV proves the universal closure of the
following disjunction with $r$ disjuncts: 
$$
[f_1(\un) = C^1_n \wedge |C^1_n| \leq \red{c}n^k \wedge (|x^1|=n\rightarrow C^1_n(x^1) = g(x^1))]\ \vee\ 
$$
$$
[f_2(\un, x^1) = C^2_n \wedge |C^2_n| \leq \red{c}n^k \wedge (|x^2|=n\rightarrow C^2_n(x^2) = g(x^2))]\ \vee\ 
$$
$$
\ldots \, \vee\  
[f_r(\un, x^1, \dots, x^{r-1}) = C^r_n \wedge |C^r_n| \leq \red{c}n^k \wedge (|x^r|=n\rightarrow 
C^r_n(x^r) = g(x^r))]\ .
$$
We shall complete the proof of the theorem by induction on $r$. The case $r=1$
and the induction step from $r-1$ to $r$ are analogous, and we describe 
only the latter.
Our induction assumption is 
that for no polynomial time functions $f_1', \dots, f_{r-1}'$ is the disjunction of the form above but with only $r-1$ disjuncts and $n$ large enough provable in PV. 

Assume without loss of generality that $k \geq 3$. By Lemma \ref{7} applied to $f \eqdef f_1$ \red{and an arbitrary but fixed $\tilde{c}_1 \geq 1$}, i.e., using the
extra hypothesis
\begin{equation}\label{6c}
\mbox{PV}\ \vdash\ \red{\mbox{UP}_{k,\tilde{c}_1}(\tilde f_1)}\ , 
\end{equation}
there is \blue{a constant $c_1 \geq 1$ and} a polynomial time algorithm $e_1$
such that \red{for $n \geq c_1$}
$$
|e_1(\un)|=n \wedge C^1_n(e_1(\un)) \neq g(e_1(\un))\ .
$$
Substitute $x^1 \eqdef e_1(\un)$ in the disjunction above. That gives \red{for large enough $n$} a \blue{valid} disjunction
of the same form (for different polynomial time functions in place of the $f_i$'s),
but with $r-1$ disjuncts:
$$
[f_2(\un, e_1(\un)) = C^2_n \wedge |C^2_n| \leq \red{c}n^k \wedge (|x^2|=n\rightarrow C^2_n(x^2) = g(x^2))]\ \vee\ 
$$
$$
\ldots \vee\  
[f_r(\un, e_1(\un),\red{x^{2}}, \dots, x^{r-1}) = C^r_n \wedge |C^r_n| \leq \red{c}n^k \wedge (|x^r|=n\rightarrow 
C^r_n(x^r) = g(x^r))]\ .
$$
This violates the induction assumption, and 
completes the induction step. 

In the proof we have used the hypothes\red{e}s that PV proves \red{$\mbox{UP}_{k,c}(g)$ for some $c \geq 1$}, \red{$\mbox{UP}_{k,\tilde{c}_1}(\tilde f_1)$ for some $\tilde{c}_1 \geq 1$}, \red{$\mbox{UP}_{k,\tilde{c}_2}$ for $\tilde f_2(1^{(n)}, e_1(1^{(n)}))$}, etc., all together $r+1$ such assumptions.
Hence one of them must fail.
This completes the proof of Theorem \ref{t:main}.
\qed

Making the informal exposition from Section \ref{s:statement} a bit more precise, observe that we do not obtain a hard language directly from a proof of $\upkconstant(g)$. This is done via a iterative process that depends on the provability of additional sentences. 

For the reader familiar with the argument in \cite[Theorem 1.1]{SW}, notice that 
we crucially used that the second application of their initial assumption does not require the uniformity condition. Roughly speaking, this would lead to the consideration of the provability in PV of a sentence expressing a uniform circuit upper bound, while here we are concerned with non-uniform circuit complexity.

\red{Finally, regarding extending Theorem \ref{t:main} to stronger theories, we remark that in Buss's theory $S^1_2$ the analogue of Theorem \ref{t:KPT} requires a number $r$ of disjuncts that may depend on $n$, and our induction on parameter $r$ could lead to (composed) functions of super-polynomial complexity.}

~\\
\red{\noindent \textbf{Acknowledgements.} We thank
Emil Je\v r\'abek for comments on an initial draft of the paper, which led us to a more robust statement of the result.}

%
%
%
%
%

%

\end{document}